\documentclass[12pt]{article}

       \usepackage{amsmath,amsthm}
       \usepackage{amssymb,amsbsy}
       \usepackage{amsfonts}

\usepackage{times}
\usepackage{cite}

\begin{document}

\author{
{\large {\bf S. Duplij${}^{\left. 1\right) }$} and  {\bf
\framebox{W. Marcinek} ${}^{\left. 2 \right) }$}}\\
 ${}^{\left. 1\right) }$ {\it Department of Physics and Technology,}\\Kharkov
National University, Kharkov 61077, Ukraine\\\emph{E-mail}:
Steven.A.Duplij@univer.kharkov.ua\\\emph{Internet}:
http://www.math.uni-mannheim.de/\~{}duplij\\
 ${}^{\left. 2\right) }${\it Institute of Theoretical Physics,}\\University
of Wroc\l aw, Pl. Maxa Borna 9, 50-204 Wroc{\l}aw, Poland
}

\title{{\bf BRAID SEMISTATISTICS AND DOUBLY REGULAR $R$-MATRIX}
\thanks{This is our unfortunately last common article (unfinished,
but contained the
main ideas)
because second author unexpectedly passed away on June 9, 2003.}
\thanks{Published in the W. Marcinek's Memorial volume of the
Kharkov University Journal Vestnik KNU,
ser. "Nuclei, particles and fields".
- 2003. - V. 601. - N 2(22). - p. 75--80.}
}

\date{June 6, 2003}

\maketitle

\noindent

\begin{abstract}
We introduce ``noninvertible\textquotedblright\ generalization of
statistics - semistatistics replacing condition when double exchanging gives identity to
``regularity'' condition. Then in categorical language we correspondingly generalize
braidings and the quantum Yang-Baxter equation. We define the
doubly regular $R$-matrix and introduce obstructed regular bialgebras.
\smallskip


\end{abstract}
\newpage

Particle systems endowed with generalized statistics and its quantizations
have been studied from different points of view (for review see e.g.
\cite{marc0,marc1}). The color statistics have been considered in \cite{marc}
(and refs. therein), and the category for color statistics has been described
in details in \cite{marc2}. The statistics in low dimensional spaces is based
on the notion of the braid group \cite{maj,maj2} (see also \cite{gre/ser} for
its acyclic extension). The construction for the category corresponding to a
given triangular solution of the quantum Yang-Baxter equation has been given
by Lyubashenko \cite{lyu3}.The statistics corresponding for arbitrary
triangular solution of the quantum Yang-Baxter equation called $S$-statistics
has been discussed by Gurevich \cite{gur1}. The mathematical formalism for the
description of particle system with $S$-statistics is based on the theory of
the tensor (monoidal) symmetric categories of MacLane \cite{maclane}. The
mathematical formalism related to an arbitrary braid statistics has been
developed by Majid \cite{maj1}. Such formalism is based on the concept of
quasitensor (braided monoidal) categories which has been introduced by Joyal
and Street \cite{joy/str}.

The previous generalizations are \textquotedblleft
invertible\textquotedblright\ in the following sense: having the two-particle
exchange process $12\rightarrow21$ (which in the simplest case usually yields
the phase factor $\pm1$ or general anyonic factor \cite{majid}), then double
exchanging gives identity $12\rightarrow12$. Here we weaken this requirement
by moving to nearest \textquotedblleft noninvertible\textquotedblright%
\ generalization of statistics -- \textquotedblleft
regularity\textquotedblright\ as follows (symbolically)
\begin{align}
12\overset{a}{\rightarrow}21\overset{b}{\rightarrow}12  &  =12\overset
{\operatorname*{id}\ }{\rightarrow}12\text{ \textquotedblleft
invertibility\textquotedblright,}\label{i}\\
12\overset{a}{\rightarrow}21\overset{b}{\rightarrow}12\overset{a}{\rightarrow
}21  &  =12\overset{a}{\rightarrow}21\text{ \textquotedblleft left
regularity\textquotedblright,}\label{r1}\\
21\overset{b}{\rightarrow}12\overset{a}{\rightarrow}21\overset{b}{\rightarrow
}12  &  =21\overset{b}{\rightarrow}12\text{ \textquotedblleft right
regularity\textquotedblright.} \label{r2}%
\end{align}

In this consideration we can treat usual statistics as \textit{one} morphism
$a$, in other words, the representation of the morphism $a$ (because $b$ can
be found from the \textquotedblleft invertibility\textquotedblright\ condition
(\ref{i}) which is $a\circ b=\operatorname*{id}$\ symbolically) by various
phase factors or elements of $R$-matrix. Here we introduce the more abstract
concept of \textquotedblleft semistatistics\textquotedblright\ as a
\textit{pair} of exchanging morphisms $a$ and $b$ satisfying the
\textquotedblleft regularity\textquotedblright\ conditions (\ref{r1}%
)--(\ref{r2}) (symbolically $a\circ b\circ a=a$, $b\circ a\circ b=$ $b$). The
general regularization procedure for different systems was previously studied
in \cite{duplij,dup/mar5,dup/mar7,dup/mar3}.

We also introduce the notion of braid semistatistics and corresponding
generalization of the quantum Yang-Baxter equation.

\section{Braid semistatistics and regular Yang-Baxter equation}

Let $\mathfrak{C}$ be a directed graph with objects $\mathrm{Ob}\mathfrak{C}$
and arrows $\mathrm{Mor}\mathfrak{C}$ \cite{mitchell,maclane1}. An
$N$\textit{-regular cocycle} $(X_{1},X_{2}\ldots,f_{1},f_{2}\ldots)$ in
$\mathfrak{C}$, $N$\ $=1,2,\ldots$, is a sequence of arrows%
\begin{equation}%
\begin{array}
[c]{c}%
X_{1}\overset{f_{1}}{\longrightarrow}X_{2}\overset{f_{2}}{\longrightarrow
}\cdots\overset{f_{N-1}}{\longrightarrow}X_{N}\overset{f_{N}}{\longrightarrow
}X_{1},
\end{array}
\end{equation}
such that
\begin{equation}%
\begin{array}
[c]{c}%
f_{1}\circ f_{N}\circ\cdots\circ f_{2}\circ f_{1}=f_{1},\\
f_{2}\circ f_{1}\circ\cdots\circ f_{3}\circ f_{2}=f_{2},\\
\\
f_{N}\circ f_{N-1}\circ\cdots\circ f_{1}\circ f_{N}=f_{N}.
\end{array}
\end{equation}

We define $N$\ \textit{obstructors} by%
\begin{equation}%
\begin{array}
[c]{c}%
e_{X_{1}}^{\left(  N\right)  }:=f_{N}\circ\cdots\circ f_{2}\circ f_{1}%
\in\operatorname*{End}(X_{1}),\\
e_{X_{2}}^{\left(  N\right)  }:=f_{1}\circ\cdots\circ f_{3}\circ f_{2}%
\in\operatorname*{End}(X_{2}),\\
\vdots\\
e_{X_{N}}^{\left(  N\right)  }:=f_{N-1}\circ\cdots\circ f_{1}\circ f_{N}%
\in\operatorname*{End}(X_{N}).
\end{array}
\label{eee}%
\end{equation}

The correspondence $e_{X}^{(N)}:X_{n}\in\mathrm{Ob}\mathfrak{C}\mapsto
e_{X_{n}}^{(N)}\in\operatorname*{End}(X_{n})$, $n$\ $=1,2,\ldots,N$\ , is
called an $N$\ -regular cocycle \textit{obstruction} on
$(X_{1},X_{2},\ldots,X_{N}|f_{1},f_{2},\ldots,f_{N})$ in $\mathfrak{C}$.

Let $\mathfrak{M}$ be a monoidal category \cite{maclane1,joy/str} which
abstractly defines the braid statistics. An $N$\ -regular \textit{obstructed
monoidal category} $\mathfrak{M}_{obtr}^{\left(  N\right)  }$ can be defined
as usual, but instead of the identity $\operatorname*{id}$\ $_{X}%
\otimes\operatorname*{id}$\ $_{Y}=\operatorname*{id}$\ $_{X\otimes Y}$ we have
an obstruction structure $e_{X}^{(N)}=\{e_{X_{n}}^{\left(  N\right)  }%
\in\operatorname*{End}(X_{n});N$\ $=1,2,...\}$ satisfying the condition
\begin{equation}
e_{X_{n}\otimes Y_{n}}^{(N)}=e_{X_{n}}^{(N)}\otimes e_{Y_{n}}^{(N)}
\label{mul}%
\end{equation}
for every two $N$\ -regular cocycles $(X_{1},X_{2},\ldots,X_{N}|f_{1}%
,f_{2},\ldots,f_{N})$ and $(Y_{1},Y_{2},\ldots$, $ Y_{N}|g_{1},g_{2},\ldots
,g_{N})$.

In a monoidal category $\mathfrak{M}$ for any two objects $X,Y\in
\mathrm{ob}\mathfrak{M}$ and the product $X\otimes Y$ one can define a natural
isomorphism (\textquotedblleft braiding\textquotedblright\ \cite{joy/str}) by
$\mathrm{B}_{X,Y}:X\otimes Y\rightarrow Y\otimes X$ satisfying the
\textit{symmetry condition} (\textquotedblleft invertibility\textquotedblright)%

\begin{equation}
\mathrm{B}_{Y,X}\circ\mathrm{B}_{X,Y}=\operatorname*{id}{}_{X\otimes Y}
\label{bbi}%
\end{equation}
which formally defines $\mathrm{B}_{Y,X}=\mathrm{B}_{X,Y}^{-1}:Y\otimes
X\rightarrow X\otimes Y$. The simplest type of braiding is the usual
transposition $\tau_{X,Y}\left(  x\otimes y\right)  =y\otimes x$, where $x\in
X$, $y\in Y$. Nonsymmetric braidings in context of the noncommutative geometry
were considered in \cite{marc3,maj3} (see also \cite{gre/ser}). In the
obstructed monoidal category $\mathfrak{M}_{obstr}^{\left(  N\right)  }$ we
introduce a \textquotedblleft regular\textquotedblright\ extension of the
braidings as follows. Let $(X_{1},X_{2},\ldots,X_{N}|f_{1},f_{2},\ldots
,f_{N})$ and $(Y_{1},Y_{2},\ldots,Y_{N}|g_{1},g_{2},\ldots,g_{N})$ are regular
cocycles and $e_{X_{n}}^{(N)}$, $e_{Y_{n}}^{(N)}$ are corresponding
obstructors, then we have two sets of monoidal products of $N$-regular
cocycles $X_{1}\otimes Y_{1}$, $X_{2}\otimes Y_{2}$, \ldots$X_{N}\otimes
Y_{N}$, $f_{1}\otimes g_{1}$, $f_{2}\otimes g_{2}$, \ldots$f_{N}\otimes g_{N}%
$, and $Y_{1}\otimes X_{1}$, $Y_{2}\otimes X_{2}$, \ldots$Y_{N}\otimes X_{N}$,
$g_{1}\otimes f_{1}$, $g_{2}\otimes f_{2}$, \ldots$g_{N}\otimes f_{N}$, and
the obstructors satisfy $e_{X_{n}}^{(N)}\otimes e_{Y_{n}}^{(N)}=e_{X_{n}%
\otimes Y_{n}}^{(N)}$.

An $N$-\textit{regular (\textquotedblleft vector\textquotedblright) braiding}
$\mathrm{\vec{B}}^{\left(  N\right)  }$ is a set of (\textquotedblleft%
$n$-component\textquotedblright) maps%
\[
X_{n}\otimes Y_{n}\overset{B_{X_{n}\otimes Y_{n}}^{\left(  N\right)  ,n}%
}{\rightarrow}Y_{n}\otimes X_{n}%
\]
such that the following diagram
\[%
\begin{array}
[c]{ccccccc}%
X_{1}\otimes Y_{1} & \overset{f_{1}\otimes g_{1}}{\rightarrow} & X_{2}\otimes
Y_{2} & \overset{f_{2}\otimes g_{2}}{\rightarrow} & \ldots & \rightarrow &
X_{N}\otimes Y_{N}\\
\mathrm{B}_{X_{n}\otimes Y_{n}}^{\left(  N\right)  ,n}\downarrow &  &
\mathrm{B}_{X_{n}\otimes Y_{n}}^{\left(  N\right)  ,n}\downarrow &  &  &  &
\mathrm{B}_{X_{n}\otimes Y_{n}}^{\left(  N\right)  ,n}\downarrow\\
Y_{1}\otimes X_{1} & \overset{g_{1}\otimes f_{1}}{\rightarrow} & Y_{2}\otimes
X_{2} & \overset{g_{2}\otimes f_{2}}{\rightarrow} & \ldots & \rightarrow &
Y_{N}\otimes X_{N}%
\end{array}
\]
is commutative. Instead of the symmetry condition (\ref{bbi}) we introduce the
\textit{generalized (1-star) inverse} $N$-regular braiding $\mathrm{\vec{B}%
}^{\ast\left(  N\right)  }$ with components satisfying%

\begin{equation}
\mathrm{B}_{X_{n}\otimes Y_{n}}^{\left(  N\right)  ,n}\circ\mathrm{B}%
_{X_{n}\otimes Y_{n}}^{\ast\left(  N\right)  ,n}\circ\mathrm{B}_{X_{n}\otimes
Y_{n}}^{\left(  N\right)  ,n}=\mathrm{B}_{X_{n}\otimes Y_{n}}^{\left(
N\right)  ,n}, \label{bbb}%
\end{equation}
where in general $\mathrm{B}_{X_{n}\otimes Y_{n}}^{\ast\left(  N\right)
,n}\neq\mathrm{B}_{X_{n}\otimes Y_{n}}^{\left(  N\right)  ,n,-1}$. We call
such a category a \textquotedblleft regular\textquotedblright\ category
\cite{dup/mar5,dup/mar7} to distinct from symmetric and \textquotedblleft
braided\textquotedblright\ categories \cite{maclane1,joy/str}.

The prebraiding relations in a symmetric monoidal category are defined as
\cite{marc1,maj2,joy/str}
\begin{align}
\mathrm{B}_{X\otimes Y,Z}  &  =\mathbf{B}_{X,Z,Y}^{\mathsf{R}}\circ
\mathbf{B}_{X,Y,Z}^{\mathsf{L}},\label{bb1}\\
\mathrm{B}_{Z,X\otimes Y}  &  =\mathbf{B}_{X,Z,Y}^{\mathsf{L}}\circ
\mathbf{B}_{X,Y,Z}^{\mathsf{R}},\label{bb2}\\
\mathbf{B}_{X,Y,Z}^{\mathsf{L}}  &  =\operatorname*{id}{}_{X}\otimes
\mathrm{B}_{Y,Z},\label{bb3}\\
\mathbf{B}_{X,Y,Z}^{\mathsf{R}}  &  =\mathrm{B}_{X,Y}\otimes\operatorname*{id}%
{}_{Z}, \label{bb4}%
\end{align}
and prebraidings $\mathrm{B}_{X\otimes Y,Z}$ and $\mathrm{B}_{Z,X\otimes Y}$
satisfy (for symmetric case) the \textquotedblleft
invertibility\textquotedblright\ property
\[
\mathrm{B}_{X\otimes Y,Z}^{-1}\circ\mathrm{B}_{X\otimes Y,Z}%
=\operatorname*{id}{}_{X\otimes Y\otimes Z},
\]
where $\mathrm{B}_{X\otimes Y,Z}^{-1}=\mathrm{B}_{Z,X\otimes Y}$. In this
notations the standard \textquotedblleft invertible\textquotedblright\ quantum
Yang-Baxter equation takes the form \cite{maj2,maj3}%

\begin{equation}
\mathbf{B}_{Y,Z,X}^{\mathsf{R}}\circ\mathbf{B}_{Y,X,Z}^{\mathsf{L}}%
\circ\mathbf{B}_{X,Y,Z}^{\mathsf{R}}=\mathbf{B}_{Z,X,Y}^{\mathsf{L}}%
\circ\mathbf{B}_{X,Z,Y}^{\mathsf{R}}\circ\mathbf{B}_{X,Y,Z}^{\mathsf{L}}.
\label{yb0}%
\end{equation}

For \textquotedblleft noninvertible\textquotedblright\ braidings satisfying
regularity (\ref{bbb}) in search of the analogs of the definitions
(\ref{bb3})--(\ref{bb4}) it is naturally to exploit the obstructors $e_{X_{n}%
}^{\left(  N\right)  }$ instead of identity $\operatorname*{id}{}_{X_{n}}$
($n=1\ldots N$) which were introduced in \cite{dup/mar,dup/mar1}. They are
defined as self-mappings $e_{X_{n}}^{\left(  N\right)  }:X_{n}\rightarrow
X_{n}$ satisfying closure conditions%

\begin{align}
e_{X_{n}}^{\left(  1\right)  }  &  =\operatorname*{id}{}_{X_{n}},\label{e1}\\
e_{X_{n}}^{\left(  2\right)  }  &  =g\circ f,\label{e2}\\
e_{X_{n}}^{\left(  3\right)  }  &  =h\circ g\circ f,\label{e3}\\
&  \ldots\nonumber
\end{align}
where $g,h...$ are some morphisms (see \cite{dup/mar1} for details). Then
using the following triple maps
\begin{align*}
\mathbf{T}_{X_{n},Y_{n},Z_{n}}^{\left(  N\right)  ,n~\mathsf{L}}  &
:X_{n}\otimes Y_{n}\otimes Z_{n}\rightarrow X_{n}\otimes Z_{n}\otimes Y_{n},\\
\mathbf{T}_{X_{n},Y_{n},Z_{n}}^{\left(  N\right)  ,n~\mathsf{R}}  &
:X_{n}\otimes Y_{n}\otimes Z_{n}\rightarrow Y_{n}\otimes X_{n}\otimes Z_{n}%
\end{align*}
defined similarly to (\ref{bb3})--(\ref{bb4})%

\begin{align}
\mathbf{T}_{X_{n},Y_{n},Z_{n}}^{\left(  N\right)  ,n~\mathsf{L}}  &
=e_{X_{n}}^{\left(  N\right)  }\otimes\mathrm{B}_{Y_{n},Z_{n}}^{\left(
N\right)  ,n},\label{eb1}\\
\mathbf{T}_{X_{n},Y_{n},Z_{n}}^{\left(  N\right)  ,n~\mathsf{R}}  &
=\mathrm{B}_{X_{n},Y_{n}}^{\left(  N\right)  ,n}\otimes e_{Z_{n}}^{\left(
N\right)  }, \label{eb2}%
\end{align}
we weaken prebraiding construction (\ref{bb1})--(\ref{bb2}) in the following
way
\begin{align}
\mathrm{P}_{X_{n}\otimes Y_{n},Z_{n}}^{\left(  N\right)  ,n}  &
=\mathbf{T}_{X_{n},Z_{n},Y_{n}}^{\left(  N\right)  ,n~\mathsf{R}}%
\circ\mathbf{T}_{X_{n},Y_{n},Z_{n}}^{\left(  N\right)  ,n~\mathsf{L}%
},\label{rb1}\\
\mathrm{P}_{Z_{n},X_{n}\otimes Y_{n}}^{\left(  N\right)  ,n}  &
=\mathbf{T}_{X_{n},Z_{n},Y_{n}}^{\left(  N\right)  ,n~\mathsf{L}}%
\circ\mathbf{T}_{X_{n},Y_{n},Z_{n}}^{\left(  N\right)  ~\mathsf{R}}.
\label{rb2}%
\end{align}
Thus the corresponding \textquotedblleft noninvertible\textquotedblright%
\ analog of the Yang-Baxter equation (\ref{rb2}) is the set of
\textquotedblleft component\textquotedblright equations%
\begin{equation}
\mathbf{T}_{Y_{n},Z_{n},X_{n}}^{\left(  N\right)  ,n~\mathsf{R}}%
\circ\mathbf{T}_{Y_{n},X_{n},Z_{n}}^{\left(  N\right)  ,n~\mathsf{L}}%
\circ\mathbf{T}_{X_{n},Y_{n},Z_{n}}^{\left(  N\right)  ,n~\mathsf{R}%
}=\mathbf{T}_{Z_{n},X_{n},Y_{n}}^{\left(  N\right)  ,n~\mathsf{L}}%
\circ\mathbf{T}_{X_{n},Z_{n},Y_{n}}^{\left(  N\right)  ,n~\mathsf{R}}%
\circ\mathbf{T}_{X_{n},Y_{n},Z_{n}}^{\left(  N\right)  ,n~\mathsf{L}}.
\label{yb}%
\end{equation}

Its solutions can be found by application of the semigroup methods (see e.g.
\cite{fangli2,fangli3}). Let us construct \textquotedblleft braidings
tower\textquotedblright\ of $k$-star regular braidings, and for 1-star regular
braidings we have%

\begin{align}
\mathrm{P}_{X_{n}\otimes Y_{n},Z_{n}}^{\left(  N\right)  ,n}\circ\overset
{\ast}{\mathrm{P}}_{X_{n}\otimes Y_{n},Z_{n}}^{\left(  N\right)  ,n}%
\circ\mathrm{P}_{X_{n}\otimes Y_{n},Z_{n}}^{\left(  N\right)  ,n}  &
=\mathrm{P}_{X_{n}\otimes Y_{n},Z_{n}}^{\left(  N\right)  ,n}\label{bbbb}\\
\overset{\ast}{\mathrm{P}}_{X_{n}\otimes Y_{n},Z_{n}}^{\left(  N\right)
,n}\circ\mathrm{P}_{X_{n}\otimes Y_{n},Z_{n}}^{\left(  N\right)  ,n}%
\circ\overset{\ast}{\mathrm{P}}_{X_{n}\otimes Y_{n},Z_{n}}^{\left(  N\right)
,n}  &  =\overset{\ast}{\mathrm{P}}_{X_{n}\otimes Y_{n},Z_{n}}^{\left(
N\right)  ,n},\\
\mathrm{P}_{Z_{n},X_{n}\otimes Y_{n}}^{\left(  N\right)  ,n}\circ\overset
{\ast}{\mathrm{P}}_{Z_{n},X_{n}\otimes Y_{n}}^{\left(  N\right)  ,n}%
\circ\mathrm{P}_{Z_{n},X_{n}\otimes Y_{n}}^{\left(  N\right)  ,n}  &
=\mathrm{P}_{Z_{n},X_{n}\otimes Y_{n}}^{\left(  N\right)  ,n},\\
\overset{\ast}{\mathrm{P}}_{Z_{n},X_{n}\otimes Y_{n}}^{\left(  N\right)
,n}\circ\mathrm{P}_{Z_{n},X_{n}\otimes Y_{n}}^{\left(  N\right)  ,n}%
\circ\overset{\ast}{\mathrm{P}}_{Z_{n},X_{n}\otimes Y_{n}}^{\left(  N\right)
,n}  &  =\overset{\ast}{\mathrm{P}}_{Z_{n},X_{n}\otimes Y_{n}}^{\left(
N\right)  ,n}, \label{bbbb1}%
\end{align}
where $\overset{\ast}{\mathrm{P}}_{X_{n}\otimes Y_{n},Z_{n}}^{\left(
N\right)  ,n}$is the generalized inverse (see e.g. \cite{nashed}) for
$\mathrm{P}_{X_{n}\otimes Y_{n},Z_{n}}^{\left(  N\right)  ,n}$, and in general
case $\overset{\ast}{\mathrm{P}}_{X_{n}\otimes Y_{n},Z_{n}}^{\left(  N\right)
,n}\neq\mathrm{P}_{X_{n}\otimes Y_{n},Z_{n}}^{\left(  N\right)  ,n,~-1}$. In a
similar we can define $k$-star braidings $\mathrm{P}_{X_{n}\otimes Y_{n}%
,Z_{n}}^{\left(  N\right)  ,n,\overset{k}{\overbrace{\ast\ast\ldots\ast}}}$
($K\times N$-regular morphisms, their number is $KN$), where $k=0,1,2\ldots
K-1$ \cite{dup/mar,dup/mar3}.

\section{Regular Yang-Baxter operators}

Let we have a set of regular obstructed algebras $\left(  A_{n},m_{n}%
,e_{A_{n}}^{\left(  N\right)  }\right)  $ with multiplication $m_{n}$ and
obstructor $e_{A_{n}}^{\left(  N\right)  }:A_{n}\rightarrow A_{n}$ (see
(\ref{eee})) such that the diagram%

\[%
\begin{array}
[c]{ccccccc}%
A_{1}\otimes A_{1} & \overset{f_{1}\otimes f_{1}}{\rightarrow} & A_{2}\otimes
A_{2} & \overset{f_{2}\otimes f_{2}}{\rightarrow} & \ldots & \rightarrow &
A_{N}\otimes A_{N}\\
m_{1}\downarrow &  & m_{2}\downarrow &  &  &  & m_{N}\downarrow\\
A_{1} & \overset{f_{1}}{\rightarrow} & A_{2} & \overset{f_{2}}{\rightarrow} &
\ldots & \rightarrow & A_{N}%
\end{array}
\]
is commutative, or
\begin{equation}
e_{A_{n}}^{\left(  N\right)  }\circ m_{n}=m_{n}\circ e_{A_{n}\otimes A_{n}%
}^{\left(  N\right)  }. \label{me}%
\end{equation}

We introduce $N$ Yang-Baxter operators $R_{n}^{\left(  N\right)  }%
:A_{n}\otimes A_{n}\rightarrow A_{n}\otimes A_{n}$ which commute with
obstructors
\begin{equation}
R_{n}^{\left(  N\right)  }\circ e_{A_{n}\otimes A_{n}}^{\left(  N\right)
}=e_{A_{n}\otimes A_{n}}^{\left(  N\right)  }\circ R_{n}^{\left(  N\right)  }
\label{re}%
\end{equation}
and satisfy $N$-regular analog of the Yang-Baxter equation (set of $N$
equations)%
\begin{align}
&  \left(  e_{A_{n}}^{\left(  N\right)  }\otimes R_{n}^{\left(  N\right)
}\right)  \circ\left(  R_{n}^{\left(  N\right)  }\otimes e_{A_{n}}^{\left(
N\right)  }\right)  \circ\left(  e_{A_{n}}^{\left(  N\right)  }\otimes
R_{n}^{\left(  N\right)  }\right) \label{ybr}\\
&  =\left(  R_{n}^{\left(  N\right)  }\otimes e_{A_{n}}^{\left(  N\right)
}\right)  \circ\left(  e_{A_{n}}^{\left(  N\right)  }\otimes R_{n}^{\left(
N\right)  }\right)  \circ\left(  R_{n}^{\left(  N\right)  }\otimes e_{A_{n}%
}^{\left(  N\right)  }\right)  .\nonumber
\end{align}

We define 1-star $N$-regular obstructed Yang-Baxter operator (set of $N$
operators $R_{n}^{\left(  N\right)  \ast}$) by%
\begin{align}
R_{n}^{\left(  N\right)  }\circ R_{n}^{\left(  N\right)  \ast}\circ
R_{n}^{\left(  N\right)  }  &  =R_{n}^{\left(  N\right)  },\label{rrr1}\\
R_{n}^{\left(  N\right)  \ast}\circ R_{n}^{\left(  N\right)  }\circ
R_{n}^{\left(  N\right)  \ast}  &  =R_{n}^{\left(  N\right)  \ast}.
\label{rrr2}%
\end{align}
Similarly, one can define $k$-star operators $R_{n}^{\overset{k}%
{\overbrace{\ast\ast\ldots\ast}}}$ ($K\times N$-regular Yang-Baxter operators,
their number is $KN$), where $n=1,2\ldots N;$ $k=0,1,2\ldots K-1$
\cite{dup/mar,dup/mar1,dup/mar2}.

\section{Bialgebras and universal $R$-matrix}

An obstructed (see \cite{dup/mar5}) $N$-regular bialgebra can be defined as a
set of $N$ bialgebras $\left(  H_{n},m_{n},\Delta_{n},e_{H_{n}}^{\left(
N\right)  }\right)  $, where $H_{n}$, $\left(  n=1\ldots N\right)  $ are
linear vector spaces over $\mathbb{K}$ with multiplications $m_{n}%
:H_{n}\otimes H_{n}\rightarrow H_{n}$ and comultiplications $\Delta_{n}%
:H_{n}\rightarrow H_{n}\otimes H_{n}$, but instead of identity map we have now
$N$ obstructors $e_{H_{n}}^{\left(  N\right)  }:H_{n}\rightarrow H_{n}$
(analogies of mappings (\ref{eee})) satisfying the consistency conditions
\begin{equation}
e_{H_{n}}^{\left(  N\right)  }\circ m_{n}=m_{n}\circ e_{H_{n}\otimes H_{n}%
}^{\left(  N\right)  },\ \ \ \ \ \ \ \Delta_{n}\circ e_{H_{n}}^{\left(
N\right)  }=e_{H_{n}\otimes H_{n}}^{\left(  N\right)  }\circ\Delta_{n}.
\label{d}%
\end{equation}

The associativity and coassociativity now have the form%
\[
m_{n}\circ\left(  m_{n}\otimes e_{H_{n}}^{\left(  N\right)  }\right)
=m_{n}\circ\left(  e_{H_{n}}^{\left(  N\right)  }\otimes m_{n}\right)
,\ \ \ \ \ \left(  \Delta_{n}\otimes e_{H_{n}}^{\left(  N\right)  }\right)
\circ\Delta_{n}=\left(  e_{H_{n}}^{\left(  N\right)  }\otimes\Delta
_{n}\right)  \circ\Delta_{n}.
\]

The Yang-Baxter operators $R_{n}^{\left(  N\right)  }:H_{n}\otimes
H_{n}\rightarrow H_{n}\otimes H_{n}$ also satisfy the additional consistency
conditions (analogy of (\ref{re}))%
\[
e_{H_{n}\otimes H_{n}}^{\left(  N\right)  }\circ R_{n}^{\left(  N\right)
}=R_{n}^{\left(  N\right)  }\circ e_{H_{n}\otimes H_{n}}^{\left(  N\right)  }%
\]
and the set of $N$ Yang-Baxter equations of type (\ref{ybr}), as follows%
\begin{align}
&\left(  e_{H_{n}}^{\left(  N\right)  }\otimes R_{n}^{\left(  N\right)
}\right)  \circ\left(  R_{n}^{\left(  N\right)  }\otimes e_{H_{n}}^{\left(
N\right)  }\right)  \circ\left(  e_{H_{n}}^{\left(  N\right)  }\otimes
R_{n}^{\left(  N\right)  }\right)  \\&=\left(  R_{n}^{\left(  N\right)  }\otimes
e_{H_{n}}^{\left(  N\right)  }\right)  \circ\left(  e_{H_{n}}^{\left(
N\right)  }\otimes R_{n}^{\left(  N\right)  }\right)  \circ\left(
R_{n}^{\left(  N\right)  }\otimes e_{H_{n}}^{\left(  N\right)  }\right)
,\nonumber
\end{align}
which defines the universal obstructed $N$-regular $R$-matrix for obstructed
$N$-regular bialgebra $\left(  H_{n},m_{n},\Delta_{n},e_{H_{n}}^{\left(
N\right)  }\right)  $. We define 1-star universal obstructed $N$-regular
$R$-matrix by%
\begin{equation}
R_{n}^{\left(  N\right)  }\circ R_{n}^{\left(  N\right)  \ast}\circ
R_{n}^{\left(  N\right)  }=R_{n}^{\left(  N\right)  },\ \ \ \ \ R_{n}^{\left(
N\right)  \ast}\circ R_{n}^{\left(  N\right)  }\circ R_{n}^{\left(  N\right)
\ast}=R_{n}^{\left(  N\right)  \ast}.
\end{equation}
As above one can define $k$-star Yang-Baxter operators $R_{n}^{\overset
{k}{\overbrace{\ast\ast\ldots\ast}}}$ (set of $KN$ operators) $n=1,2\ldots N;$
$k=0,1,2\ldots K-1$ \cite{dup/mar,dup/mar1}. Then the convolution product can
be defined (in \textquotedblleft components\textquotedblright) as%
\begin{equation}
s\star_{n}t:=m_{n}\circ(s\otimes t)\circ\triangle_{n}, \label{conv}%
\end{equation}
where $s,t\in\hom_{m_{n}}(H_{n},H_{n})$.

Let $A$ be an $N$-regular obstructed algebra with $N$ obstructors
$e_{A}^{\left(  N\right)  }$ and multiplication $m$, and $R^{\left(  N\right)
}$ be an $N$-regular Yang-Baxter operator on $A$, then the algebra $A$ with
the multiplication $m_{R}=m\circ R^{\left(  N\right)  }$ is also an
$N$-regular obstructed algebra. Indeed, from definition (\ref{me}) we have
$e_{A}^{\left(  N\right)  }\circ m=m\circ e_{A\otimes A}^{\left(  N\right)  }%
$, and then from (\ref{re}) we obtain $m_{R}\circ e_{A\otimes A}^{\left(
N\right)  }=m\circ R^{\left(  N\right)  }\circ e_{A\otimes A}^{\left(
N\right)  }=m\circ e_{A\otimes A}^{\left(  N\right)  }\circ R^{\left(
N\right)  }=e_{A}^{\left(  N\right)  }\circ m\circ R^{\left(  N\right)
}=e_{A}^{\left(  N\right)  }\circ m_{R}$.

Let $C$ be an $N$-regular obstructed coalgebra with $N$ obstructors
$e_{C}^{\left(  N\right)  }$ and comultiplication $\Delta$, and $R^{\left(
N\right)  }$ be an $N$-regular Yang-Baxter operator on $C$, then the algebra
$C$ with the comultiplication $\Delta_{R}=R^{\left(  N\right)  }\circ\Delta$
is also an $N$-regular obstructed coalgebra. Indeed, from definition (\ref{d})
we have $\Delta\circ e_{A}^{\left(  N\right)  }=e_{A\otimes A}^{\left(
N\right)  }\circ\Delta$, and then from (\ref{re}) we obtain $\Delta_{R}\circ
e_{A}^{\left(  N\right)  }=R^{\left(  N\right)  }\circ\Delta\circ
e_{A}^{\left(  N\right)  }=R^{\left(  N\right)  }\circ e_{A\otimes A}^{\left(
N\right)  }\circ\Delta=e_{A\otimes A}^{\left(  N\right)  }\circ R^{\left(
N\right)  }\circ\Delta=e_{A\otimes A}^{\left(  N\right)  }\circ\Delta_{R}$.

\section{Doubly regular Hopf algebras}

Usual antipode is defined as inverse to the identity under convolution, if and
only if there exist unit and counit for a bialgebra \cite{abe,sweedler}. Since
we do not require existence of unit and counit in obstructed bialgebras, we
have to define some more general analog of antipode. The Von Neumann regular
antipode for weal Hopf algebras was considered in
\cite{dup/li1,dup/li2,dup/li3} (\textquotedblleft non-unital\textquotedblright%
/\textquotedblleft nonsymmetric\textquotedblright\ antipodes were considered
in \cite{boh/nil/szl}). By analogy we can introduce the obstructed $N$-regular
antipode (set of $N$ antipodes) for every bialgebra $\left(  H_{n}%
,m_{n},\Delta_{n},e_{H_{n}}^{\left(  N\right)  }\right)  $ as a generalized
inverse for obstructor%
\begin{equation}
e_{H_{n}}^{\left(  N\right)  }\star_{n}S_{n}^{\left(  N\right)  }\star
_{n}e_{H_{n}}^{\left(  N\right)  }=e_{H_{n}}^{\left(  N\right)  }%
,\ \ \ \ S_{n}^{\left(  N\right)  }\star_{n}e_{H_{n}}^{\left(  N\right)
}\star_{n}S_{n}^{\left(  N\right)  }=S_{n}^{\left(  N\right)  }.\label{se2}%
\end{equation}

In this way we define $LN$ higher $L$-regular analogs of antipode $S_{_{n}%
}^{\overset{l}{\overbrace{\ast\ast\ldots\ast}}}$ ($l=0,1,2\ldots L-1$),
similarly to $K$-star regular quantities above. For example, in the case $l=1$
we have instead of (\ref{se2}) the following set of defining equations%
\begin{align*}
e_{H_{n}}^{\left(  N\right)  }\star_{n}S_{n}^{\left(  N\right)  }\star
_{n}S_{n}^{\left(  N\right)  \ast}\star_{n}e_{H_{n}}^{\left(  N\right)  }  &
=e_{H_{n}}^{\left(  N\right)  },\\
S_{n}^{\left(  N\right)  }\star_{n}S_{n}^{\left(  N\right)  \ast}\star
_{n}e_{H_{n}}^{\left(  N\right)  }\star_{n}S_{n}^{\left(  N\right)  }  &
=S_{n}^{\left(  N\right)  },\\
S_{n}^{\left(  N\right)  \ast}\star_{n}e_{H_{n}}^{\left(  N\right)  }\star
_{n}S_{n}^{\left(  N\right)  }\star_{n}S_{n}^{\left(  N\right)  \ast}  &
=S_{n}^{\left(  N\right)  \ast}.
\end{align*}

An obstructed $N$-regular bialgebra $\left(  H_{n},m_{n},\Delta_{n},e_{H_{n}%
}^{\left(  N\right)  }\right)  $ with $L$-regular antipode is called
obstructed $N\times L$-regular (doubly regular) Hopf algebra $$\left(
H_{n},m_{n},\Delta_{n},e_{H_{n}}^{\left(  N\right)  },S_{_{n}}^{\overset
{l}{\overbrace{\ast\ast\ldots\ast}}}\right)  ,$$ where $n=1,2\ldots N;$
$l=0,1,2\ldots L-1$.

Note, that in general, obstructed $N\times L$-regular Hopf algebras $$\left(
H_{n},m_{n},\Delta_{n},e_{H_{n}}^{\left(  N\right)  },S_{_{n}}^{\overset
{l}{\overbrace{\ast\ast\ldots\ast}}}\right)  $$ do not contain unit and/or
counit (analogously to \cite{dup/li3,boh/nil/szl}).

In the opposite case it can be possible that for each $N\times L$-regular Hopf
algebra $$\left(  H_{n},m_{n},\Delta_{n},e_{H_{n}}^{\left(  N\right)  }%
,S_{_{n}}^{\overset{l}{\overbrace{\ast\ast\ldots\ast}}}\right)  $$ there exist
unit $\eta_{n}$ and counit $\varepsilon_{n}$. If we have one antipode for
every $n$, then it should satisfy%

\[
\left(  S_{_{n}}^{\left(  N\right)  }\otimes e_{H_{n}}^{\left(  N\right)
}\right)  \circ\Delta_{n}=\left(  e_{H_{n}}^{\left(  N\right)  }\otimes
S_{_{n}}^{\left(  N\right)  }\right)  \circ\Delta_{n}=\eta_{n}\circ
\varepsilon_{n}.
\]

We call $P_{n},Q_{n}$ obstructed $N$-regular modules if for each $n$ there
exist maps $\rho_{P_{n}}:P_{n}\otimes H_{n}\rightarrow P_{n}$ and $\rho
_{Q_{n}}:H_{n}\otimes Q_{n}\rightarrow Q_{n}$, such that%
\[
e_{P_{n}}^{\left(  N\right)  }\circ\rho_{P_{n}}=\rho_{P_{n}}\circ\left(
e_{P_{n}}^{\left(  N\right)  }\otimes e_{H_{n}}^{\left(  N\right)  }\right)
,\ \ \ \ \rho_{Q_{n}}\circ e_{Q_{n}}^{\left(  N\right)  }=\left(  e_{H_{n}%
}^{\left(  N\right)  }\otimes e_{Q_{n}}^{\left(  N\right)  }\right)  \circ
\rho_{Q_{n}},
\]
where $e_{P_{n}}^{\left(  N\right)  }$ and $e_{Q_{n}}^{\left(  N\right)  }$
are obstructors for modules $P_{n}$ and $Q_{n}$ (see (\ref{eee})).

Let $R_{n}$ be the universal obstructed $N$-regular $R$-matrix on the
obstructed $N$-regular bialgebra $\left(  H_{n},m_{n},\Delta_{n},e_{H_{n}%
}^{\left(  N\right)  }\right)  $, and $P_{n}$ and $Q_{n}$ are left modules
over $H_{n}$, then there is obstructed $N$-regular braiding $B_{P_{n},Q_{n}%
}^{\left(  N\right)  }:P_{n}\otimes Q_{n}\rightarrow Q_{n}\otimes P_{n}$, such
that $B_{P_{n},Q_{n}}^{\left(  N\right)  }\left(  P_{n}\otimes Q_{n}\right)
=\tau_{P_{n},Q_{n}}\left(  R_{n}^{\left(  N\right)  }\left(  P_{n}\otimes
Q_{n}\right)  \right)  $, where $R_{n}^{\left(  N\right)  }$ is the
corresponding Yang-Baxter operator.

\section{Conclusions}

Thus, in this paper we have constructed a general categorical approach for
systems endowing \textquotedblleft noninvertible\textquotedblright%
\ (\textquotedblleft regular\textquotedblright) statistics (\ref{r1}%
)--(\ref{r2}) --- semistatistics --- using methods of \cite{marc0}%
-\cite{dup/mar3}. We introduced doubly regular prebraiding and braiding and
obtained the set of regular Yang-Baxter equations in terms of obstructors. The
doubly regular Yang-Baxter operators, bialgebras and Hopf algebras are considered.

 \end{document}